\documentclass[a4paper,11pt]{article}
\usepackage{paralist}   

\usepackage{amsmath}
\usepackage{amssymb}
\usepackage{amsfonts}
\usepackage{color}
\usepackage{multicol}
\raggedcolumns
\usepackage{verbatim}

\usepackage{graphicx}
\usepackage{subfig}
\usepackage{bm}
\usepackage{tikz}

\usepackage{lipsum}
\usepackage{ntheorem}

\renewcommand{\leq}{\leqslant}

\renewcommand{\geq}{\geqslant}
\newtheorem{theorem}{Theorem}[section]
\renewtheorem*{theorem*}{Theorem}

\newtheorem{corollary}[theorem]{Corollary}

\newtheorem{lemma}[theorem]{Lemma}
\renewtheorem*{remark*}{Remark}%[section]

\newcommand{\proof}{\textbf{Proof:\ }}
\newcommand{\pbox}{\hfill$\Box$\\}

\newcommand{\B}{\mathcal{B}}
\newcommand{\R}{\mathbb{R}}
\newcommand{\N}{\mathbb{N}}
\newcommand{\C}{\mathbb{C}}

\newcommand{\Z}{\mathbb{Z}}
\renewcommand{\S}{\mathcal{S}}

\newcommand{\D}{\mathbb{D}}
\newcommand{\I}{\mathcal{I}}
\newcommand{\F}{\mathcal{F}}

\usepackage{geometry}
\usepackage[bookmarks,bookmarksopen=false,% Klappt die Bookmarks in Acrobat aus
                  pdfauthor={Autor},%
                  pdftitle={Titel},%
                  colorlinks=true,%
                  linkcolor=blue,%
                  citecolor=blue,%
                  urlcolor=blue,%
                ]{hyperref}
                

\begin{document}

\title{Planar sampling sets %and Logvinenko-Sereda theorems 
for the short-time Fourier transform}
\author{Philippe Jaming\thanks{%
philippe.jaming@math.u-bordeaux.fr}\ \ and Michael Speckbacher\thanks{%
michael.speckbacher@u-bordeaux.fr} \vspace{0.2cm} \\ 
\emph{Univ. Bordeaux, CNRS, Bordeaux INP, IMB, UMR 5251}\\  \emph{F-33400, Talence, France}}
\date{}
\maketitle

\begin{abstract}
\noindent This paper considers the problem of restricting the short-time Fourier transform to domains of nonzero measure in the plane and studies  sampling bounds of such systems.  In particular, we give a quantitative estimate for the lower sampling bound in the case of Hermite windows and derive a sufficient condition for a large class of windows in terms of a certain planar density. On the way, we prove a Remez-type inequality for polyanalytic functions.
\end{abstract}


\section{Introduction}
In this paper, we investigate the existence and behavior of lower norm bounds for the problem of restricting the short-time Fourier transform $V_g$ to a domain in $\C$ of nonzero measure. This can also be viewed as a planar subsampling problem or a concentration problem. 

More precisely, we are looking for conditions on a measurable set $\Omega\subset \C$ to be
a \emph{sampling set} (or \emph{dominating set}) for the short-time Fourier transform $V_g$ 
in the sense that there exists a constant $C$ depending only on $g$ and $\Omega$ such that
$$
  \|V_g f\|_{L^p(\C)}\leq C\|V_gf\|_{L^p(\Omega)},\qquad \forall f\in M^p(\R),
$$
where $M^p(\R)$ denotes the \emph{modulation space} associated to $L^p(\C)$.
Moreover, we want to estimate the \emph{sampling constant} $C$ that appears in this inequality in terms of the window $g$ and geometric properties of $\Omega$.

The question of existence of such a sampling constant has been addressed in different contexts during the last decades. One of the first instances of such a problem is in the context of Fourier analysis. Here, the task is to determine sets $\Omega\subset\R$ and the constant $C$ such that 
\begin{equation}\label{eq:PW-LS}
\|f\|_{L^2(\R)}\leq C \|f\|_{L^2(\Omega)},\qquad \forall f\in L^2(\R) \mbox{ satisfying }\widehat{f}(\xi)=0,\mbox{ if }|\xi|>W/2.
\end{equation}
This question has applications in signal processing, but also in control of PDEs ({\it see e.g.} \cite{BJPS}) and local solvability of PDEs ({\it see e.g.} \cite[Theorem 10.10]{MS}).
%
%\textsc{K. Beauchard, Ph. Jaming \& K. Pravda-Starov} 
%Spectral inequality for Hermite functions and null-controllability 
%of hypoelliptic quadratic equations from thick sets. 
%{\em Arxiv:1804.04895}
%
%C. Muscalu \& W. Schlagg
%Classical and multilinear harmonic analysis, vol I
%Cambridge University texts in mathematics 137, 2013.
%
The first solution was given by the well-known Logvinenko-Sereda theorem (Panejah \cite{pa62,pa66}, Kac’nelson \cite{kac73} and Logvinenko-Sereda \cite{logse73}), and  the sampling constant has since been improved by Kovrijkine \cite{kov01} to an  essentially optimal quantitative estimates, see also \cite{rez10}.
Later, this approach was adapted to derive estimates of the sampling constant for  Bergman spaces  \cite{lue81}, functions with compactly supported Fourier-Bessel transform \cite{ghoja13}, model spaces \cite{hajake17} and finite expansions on compact manifolds \cite{orpri13}.

Since discrete sampling sets for the Paley-Wiener space are required to satisfy a certain density condition, it is not surprising that this remains true for non-discrete sampling sets. In particular, the existence of a sampling constant $C$ in \eqref{eq:PW-LS} is equivalent to the property that each interval of a given fixed size contains at least a minimum fraction of the sampling set.
More precisely, the validity of \eqref{eq:PW-LS} is equivalent to $\Omega$ being 
\emph{relatively dense}. 

Let us denote by $\D(z,R)\subset \C$ the disc of radius $R>0$ centered at $z\in\C$.  Recall that a measurable set $\Omega\subset \C$  is called 
\emph{$(\gamma,R)$-dense} if  
\begin{equation}\label{eq:gR-dense-D}
\gamma=\inf_{z\in\C}\frac{|\Omega\cap \D(z,R)|}{|\D(z,R)|}>0,
\end{equation}
and \emph{relatively dense} if there exist $\gamma,R>0$ such that $\Omega$ is $(\gamma,R)$-dense.
Janson, Peetre and Rochberg \cite{japero87} and  Ortega-Cerd{\`a} \cite{orce98}  proved that $\Omega\subset \C$ is a planar set of sampling for the Bargmann-Fock space of analytic functions if and only if $\Omega$ is relatively dense. As a direct consequence, this result settles the question of determining
planar sets of sampling for the short-time Fourier transform with Gaussian window. 
%We will also say that   $\Omega$ has density $\gamma$ at scale $R$. 
Later, Ascensi extended this characterization to a larger class of window functions that are nonzero almost everwhere and satisfy certain decay conditions %on the window and its associated reproducing kernel 
\cite[Section 6.1]{Ascensi}.

However, all  results up to this date are non-quantitative and thus do not provide estimates of the sampling constant.  
Moreover, the equivalence of planar sets of sampling and relatively dense sets cannot be extended to general window functions. Although each planar set of sampling is necessarily relatively dense \cite[Theorem 10]{Ascensi}, the opposite is not true. Take for example $g,f$ to be compactly supported. Then $V_gf$ is supported on a strip $\mathfrak{S}=[-S,S]+i\R$ in phase space. Taking $\Omega=\C\backslash \mathfrak{S}$, we see that $\Omega$ cannot be a planar set of sampling but for every $R>2S$, there exists $\gamma>0$  such that $\Omega$  is $(\gamma,R)$-dense. 

It is the main goal of this contribution to establish quantitative estimates of the sampling bounds in the case of Hermite functions and to investigate under which conditions a $(\gamma,R)$-dense set (at specific scales $R$) is also a planar set of sampling. 
 
Our main result is the following. It gives quantitative estimates for the case of Hermite windows 
(see \eqref{eq:defHerm} for the definition)
and can also be formulated in terms of true polyanalytic functions, see Corollary~\ref{cor:main}. 

\begin{theorem}\label{thm:main}
Let $1\leq p<\infty$, $\gamma,R>0$ and $h_n$ be the $n$-th Hermite function.
Then there exists $\eta=\eta(n,R)$, $\sigma=\sigma(n,R)$ and a numerical constant $C>0$
such that if $\Omega\subset \C$ is $(\gamma,R)$-dense for some scale $R>0$ and
 $f\in M^p(\R)$, then
\begin{equation}
\|V_{h_n}f\|_{L^p(\C)}\leq\eta\left(\frac{\gamma}{C}\right)^{-\sigma}\|V_{h_n}f\|_{L^p(\Omega)}.
\end{equation}
\end{theorem}
A more precise estimate of $\eta$ and $\sigma$ will be given below. The proof follows in parts  the strategy
of Kovrijkine for  the Paley-Wiener space \cite{kov01}.
The two key ingredients of this proof are  Bernstein's inequality and Remez' inequality.
In the course of proving Theorem~\ref{thm:main}, we derive a Remez-type inequality for polyanalytic functions (Theorem~\ref{prop-poly-remez}) from a result on plurisubharmonic functions \cite{br99}.
On the other hand, the lack of  Bernstein's inequality in $\C$ is overcome by a local reproducing formula for the short-time Fourier transform with Hermite functions \cite{abspe18-sieve},
 a generalization of Seip's formula for the Bargmann Fock space \cite{sei91}.

In addition, given a window function $g\in L^2(\R)$, we give a partial answer to the following problem:
\begin{center}
\emph{Does there exist $R(g)>0$ such that any $(\gamma,R)$-dense set $\Omega$  is \\ a planar set of sampling  for $R\leq R(g)$,}
\end{center}
 and give an estimate of the sampling constant in terms of $\gamma,R,$ and $g$. We show that under mild conditions on the window function, there exists a scale $R^*$ such that the above holds:
 
\begin{theorem}\label{thm:main2}
Let $g\in H^1(\R)$ be compactly supported, then there exists $R(g)$ and $\kappa=\kappa(g)$
such that if $\Omega\subset \C$ is $(\gamma,R)$-dense at scale $R\leq R(g)$ and
 $f\in L^2(\R)$, then
\begin{equation}
\|V_{g}f\|_{L^2(\C)}\leq\frac{\kappa}{\gamma^{1/2}}\|V_{g}f\|_{L^2(\Omega)}.
\end{equation}
\end{theorem}

\noindent As an application, let us mention that, if $p\in\C[X,Y]$ is a polynomial in two variables, then
for every $R>0$, there exists $\varepsilon,\gamma$ (depending only on $R$ and $p$)
such that the level set $\Omega:=\{z\in\C\,:|p(z,\overline{z})|\geq \varepsilon\}$ is
$(\gamma,R)$-dense (see e.g. \cite[Section 10.4.2]{MS}).
But then
$$
\|V_{g}f\|_{L^2(\Omega)}\leq \frac{1}{\varepsilon}\|pV_{g}f\|_{L^2(\Omega)}\leq 
\frac{1}{\varepsilon}\|pV_{g}f\|_{L^2(\C)}.
$$
Together with Theorem \ref{thm:main} or  Theorem \ref{thm:main2}, we obtain the following version
of Heisenberg's inequality for the short-time Fourier transform (see e.g. \cite{BDJ,groe1} for other versions of Heisenberg's inequality for the STFT):

\begin{corollary}
Let $g\not=0$ be either $g=h_n$ or $g\in H^1(\R)$ with compact support. Let $p$ be a polynomial of two variables.
Then there exists a constant $C=C(g,p)$ such that, for every $f\in L^2(\R)$
$$
\int_{\C}|p(z,\overline{z})V_gf(z)|^2dz\geq C\|f\|_2^2.
$$
\end{corollary}
 
%This result can also be formulated for true polyanalytic functions, see Corollary~\ref{cor:main} and we will extend it to the case (Theorem\ref{thm:multiplexed}) of  multiplexed short-time Fourier transforms, i.e.  if $\bm{f}:=(f_0,...,f_n)$ and $\bm{f}:=(h_0,...,h_n)$ then
%$$
%\bm{V_{h_n}f}(z):= \sum_{k=0}^n V_{h_k}f_k(z).
%$$

%\begin{corollary}
%Let $1\leq p<\infty$, and $F\in \bm{F}^n_p(\R)$.
%If $\Omega\subset \C$ is  is $(\gamma,R)$-dense for some scale $R>0$, then there exists $C\geq 1$ $($independent of $R,n$, and $F)$ such that 
%\begin{equation}
%\|F\|_{\mathcal{L}^p(\C)}\leq \frac{R^2}{\nu_n(R)}C^{R^2}\left(\frac{\gamma}{C}\right)^{- C \cdot \sigma(R,n)}\|F\|_{\mathcal{L}^p(\Omega)}.
%\end{equation}
%\end{corollary}

%\red{in the proof we derive a Remez type inequality on polyanalytic functions}

%\red{ This is particularly interesting for applications where one wants to measure $f$
%on a set  as small as possible while aiming at an estimate closest possible to the norm of
%$f$.  This requires a knowledge of the sampling constants depending on the size of.  For
%the Paley-Wiener space, essentially optimal quantitative estimates of these constants were
%given by Kovrijkine \cite{kov01}, see also \cite{rez10}.} 

\noindent Alternatively,  sampling bounds  can be obtained from upper bounds on the concentration problem for the complement of $\Omega$. The use of large sieve methods has been introduced for the estimate of the sampling constant in \eqref{eq:PW-LS} by Donoho-Logan \cite{dolo92}. This approach has recently been extended to the short-time Fourier transform with Hermite window, see
\cite{abspe17-sampta,abspe18-sieve} and can also be adapted to finite spherical harmonic expensions on the sphere \cite{hryspe19}. If the sets are ``thicker'' than generic dense sets, this leads to better constants than the one in this paper. Further results in this direction can be found in \cite{fega10}.
% established estimates if the set $\Omega^c$ is a thin sets at infinity.

\medskip

The paper is organized as follows. The next section is devoted to preliminaries. In Section 3, we prove a Remez inequality for polyanalytic functions. Section 4 is devoted to the proof of Theorem \ref{thm:main}
while Section 5 is devoted to the proof of Theorem \ref{thm:main2}.

\section{Preliminaries and Notation}
Throughout this paper we will write $Q_R:=\{x+i\xi\in\C:\ \max(|x|,|\xi|)\leq R/2\}$ for the square with sidelength $R$, $\D(z,R)$ for the disc in $\C$ with radius $R>0$ and center $z$, $B_\C(z,R)$ for the ball in $\C^d$, and $B_\R(z,R)$ for the restriction of $B_\C(z,R)$ to $\R^d$. Moreover,
we use the following convention for the Fourier transform
$$
\widehat{f}(\xi)=\int_\R f(t)e^{-2\pi i\xi t}dt,
$$
and define the Hermite functions by 
\begin{equation}
\label{eq:defHerm}
h_n(t)=c_n e^{\pi t^2}\left(\frac{d}{dt}\right)^n \left(e^{- 2\pi t^2}\right),
\end{equation}
where $c_n$ is chosen such that $\|h_n\|_2=1$.

\subsection{The Short-Time Fourier Transform}

Let $z=x+i\xi\in\C$. The \emph{time-frequency shift} $\pi(z)$ of a function $g:\R\rightarrow\C$ is defined as 
$$\pi(z)g(t):=M_\xi T_xg(t)=e^{2\pi i\xi t}g(t-x),$$ 
where $T_xg(t)=g(t-x)$ denotes the translation, and $M_{\xi}g(t)=e^{2\pi i\xi t}g(t)$ the modulation operator. For $f,g\in L^2(\R)$, the \emph{short-time Fourier transform} of $f$ with window $g$ defined as 
$$
V_gf(z):=\langle f,\pi(z)g\rangle_{L^2(\R)}=\int_\R f(t)\overline{g(t-x)}e^{-2\pi i\xi t}dt,
$$ 
%The short-time Fourier transform can be written as a Fourier transform in the frequency variable
%$$
%V_gf(x,\xi)=\F(f\cdot T_x)(\xi),
%$$
 is a scalar multiple of an isometry from $L^2(\R)$ to $L^2(\C)$. In particular,
$$
\int_{\C}|V_g f(z)|^2dz=\|f\|^2_2\|g\|^2_2,\qquad \forall f,g\in L^2(\R).
$$
Note that we may also define $V_gf(z):=\langle f,\pi(z)g\rangle_{L^2(\R)}$ when $g\in\S(\R)$ (the Schwartz class) and
$f\in\S'(\R)$ (a tempered distribution) and that $V_gf$ is then a locally bounded function.

As we derive lower bounds for general $L^p$-spaces, we need  to recall the definition of modulation spaces which were introduced by Feichtinger \cite{Fei-Mod}. Following \cite{groe1}, one can define the modulation space $M^p(\R)$,  $1\leq p\leq \infty$
as the space of all tempered distributions $f$ for which
$$
\|f\|_{M^p(\R)}:=\|V_{h_0}f\|_{L^p(\C)}
$$
is finite. Note that $\|V_{h_n}f\|_{L^p(\C)}$ is an equivalent norm on $M^p(\R)$ for any $n\in\N_0$.
For further reading on time-frequency analysis we refer to the standard textbooks \cite{fo89,groe1}.

\subsection{Hermite Windows and Spaces of Polyanalytic Functions}
A function $F:\C\rightarrow \C$ is called \emph{polyanalytic of order $n$} if it satisfies the higher order Cauchy-Riemann equation $(\bar\partial)^{n+1}F=0$. In that case, $F$ can be written as 
\begin{equation}\label{poly1}
F(z)=F(z,\overline{z})=\sum_{k=0}^n F_k(z)\overline{z}^k,
\end{equation}
where $F_0,\ldots,F_n:\C\rightarrow\C$ are holomorphic functions.

The \emph{true polyanalytic Bargmann transform} $\mathcal B^{n+1}$ of  a function $f\in L^2(\R)$  is defined via the short-time Fourier transform of $f$ using Hermite window $h_n$, see \cite[Section 2.2]{abgroe12}:
\begin{equation}\label{def-poly-bargmann-trafo}
\mathcal{B}^{n+1}f(z)=V_{h_n}f(\overline{z})e^{\pi (z^2-\overline{z}^2)/4} e^{\pi |z|^2/2}.
\end{equation}
In particular, $|\B^{n+1}f(z)|=|V_{h_n}f(\overline{z})|e^{\pi |z|^2/2}$.
The \emph{polyanalytic Bargmann} space $\bm{F}^n_p(\C)$ is defined as the space of all polyanalytic functions of order $n$ such that 
\begin{equation}\label{def:LpC}
\|F\|_{\mathcal{L}^p}^p:=\int_\C |F(z)|^pe^{-\pi p|z|^2/2}dz<\infty.
\end{equation}
Moreover, the \emph{true polyanalytic Bargmann} space $\F^n_p(\C)$ is the subspace of $\bm{F}^n_p(\C)$ consisting of all those functions $F$ for which there exists an analytic function $H$ such that 
$$
F(z)=\left(\frac{\pi^{n}}{n!}\right)^{\frac{1}{2}}e^{\pi|z|^2}\left(\frac{d}{dz}\right)^{n}\left[e^{-\pi|z|^2}H(z)\right].
$$
It was shown in \cite[Section 3.2 and 3.3]{abgroe12} that the images of the true polyanalytic Bargmann transform applied to the modulation spaces $M^p(\R)$ are simply the true polyanalytic Bargmann spaces $\F_p^n(\C)$, i.e. $\B^{n+1}(M^p(\R))=\F_p^{n}(\C)$, and in particular,
$\B^{n+1}(L^2(\R))=\F_2^{n}(\C)$. Moreover, the polyanalytic Bargmann spaces can be written as  the direct sum of the true polyanalytic Bergmann spaces, i.e. for $1\leq p<\infty$
$$
\bm{F}^n_p(\C)=\bigoplus_{k=0}^n \F_p^k(\C).
$$
Let $L_n$ be the $n$-th Laguerre polynomial given by the closed form $L_n(t)=\sum_{k=0}^n\binom{n}{k}\frac{(-1)^k}{k!}t^k$.
 In  \cite[Theorem 1]{abspe18-sieve} the following \emph{local reproducing formula} is shown 
to hold for every $f\in\S'(\R)$
\begin{equation}\label{eq:local-repr}
V_{h_n}f(z)=\nu_n(R)^{-1}\int_{\D(z,R)} V_{h_n}f(w)\langle \pi(w)h_n,\pi(z)h_n\rangle dw,
\end{equation}
where 
\begin{equation}\label{eq:local-const}
\nu_{n}(R):=\int_0^{\pi R^2} L_n(t)^2e^{-t}dt.
\end{equation}
For the case $n=0$, i.e. the case of the Gaussian window, this result can be deduced from Seip's local reproducing formula for the Bargmann-Fock space \cite{sei91}.

\subsection{Maximum Modulus Principle for Polyanalytic Functions}

For our proof of the Remez-type inequality for polyanalytic functions (Theorem~\ref{prop-poly-remez}), we need the maximum modulus principle for polyanalytic functions, see Balk \cite[Theorem 1.5]{balk97}. 

\begin{lemma}[Balk]\label{max-mod-principle}
If $F$ is a polyanalytic function of order $n$ in a disc $\D(0,\lambda R)$ for some $\lambda>1$ and $M:=\sup_{\D(0,\lambda R)}|F|$, then there exists $D_n=D_n(\lambda)>0$, only depending on $\lambda$, such that
$$
\sup_{z\in \D(0,R)}|F_k(z)z^k|\leq D_n M,\qquad k=0,\ldots,n,
$$
where $D_n$ is given by 
\begin{equation}\label{eq:def-Dn}
D_n:=\left(\frac{2\lambda}{\lambda-1}\right)^{n+2}(n+2)^{(n+2)^2}.
\end{equation}
\end{lemma} 

\noindent As the proof by Balk leaves out technical details and  does not reveal the dependence of the constants on the order $n$, we precise the arguments in the following.

A polyanalytic function $F$ of order $n$ is called \emph{reduced} if it can be written as
$$
F(z)=\sum_{k=0}^n H_k(z)|z|^{2k},
$$
where $H_k$ is a holomorphic function.
If $F$ is a reduced polyanalytic function it satifies a Cauchy-type formula \cite[Section 1.3, (11)]{balk97}.
\begin{lemma}
Let $F$ be a reduced polyanalytic function in $\D(0,R)$, $0<R_0<R_1<\ldots<R_n<R$, and let $\Gamma_k:=\{z:\ |z|=R_k\}$. For every $z\in \D(0,R_0)$, $F$ satisfies
\begin{equation}\label{eq:cauchy-formula}
F(z)=\frac{1}{2\pi i}\sum_{k=0}^n P_k(|z|^2)\int_{\Gamma_k}\frac{F(t)}{t-z}dt,
\end{equation}
where $P_k$ is a polynomial given by
\begin{equation*}
P_k(t):=\prod_{j\neq k}\frac{R_j^2-t}{R_j^2-R_k^2}.
\end{equation*}
\end{lemma} 

\begin{lemma}\label{estimate-F_k}
If $F$ is a polyanalytic function of order $n$ in $\D(0,\lambda R)$, for some $\lambda>1$, and $M=\sup_{\D(0,\lambda R)}|F|$, then for every $z\in D(0,(1+\lambda)R/2)$ it holds
\begin{equation}\label{eq:est-F_n}
|z^nF_n(z)|\leq M\left(\frac{2\lambda(n+2)}{\lambda-1}\right)^{n+2},
\end{equation}
and 
\begin{equation}\label{eq:est-F-minus}
|F(z)-\overline{z}^n F_n(z)|\leq  2M\left(\frac{2\lambda(n+2)}{\lambda-1}\right)^{n+2}.
\end{equation}
\end{lemma}
\proof First, observe that $z^nF(z)$ is a reduced polyanalytic function satisfying $\bar\partial^n (z^n F(z))$ $=z^n\bar \partial^n  F(z)=n!z^n F_n(z)$.
  If we choose $R_k=\frac{R}{2}(1+\lambda+(\lambda-1)\frac{k+1}{n+2})$, then $(1+\lambda)R/2<R_0<R_{1}<\ldots<R_{n}<\lambda R$ and 
$$
|R_j^2-R_k^2|=|(R_j+R_k)(R_j-R_k)|\geq 2\cdot \frac{(1+\lambda)R}{2}\cdot\frac{(\lambda-1)R}{2(n+2)}\geq\frac{\lambda(\lambda-1)R^2}{2(n+2)}.
$$  
As $P_k$ is a polynomial of degree $n$ with leading coefficient $(-1)^n/\prod_{j\neq k} (R_j^2-R_k^2)$, it follows that 
$$
P_k^{(n)}(z)=\frac{(-1)^n n!}{\prod_{j\neq k} (R_j^2-R_k^2)}.
$$
and consequently that $|P_k^{(n)}|\leq n!\left( \frac{2(n+2)}{\lambda(\lambda-1)R^{2}}\right)^n.$
  By \eqref{eq:cauchy-formula}, we may thus write
\begin{align*}
  |z^n F_n(z)|&=\frac{1}{n!}|\bar\partial^n (z^n F(z))|\leq \frac{1}{2\pi n!}\sum_{k=0}^n \big|\bar\partial^n P_k(|z|^2)\big|\int_{\Gamma_k}\frac{|t^nF(t)|}{|t-z|}dt
  \\
  &\leq \frac{1}{2\pi n!}\sum_{k=0}^n |z|^n \big|P_k^{(n)}(|z|^2)\big|\int_{\Gamma_k}\frac{|t^nF(t)|}{|t-z|}dt
  \\
  &\leq \frac{1}{2\pi }\sum_{k=0}^n (\lambda R)^n \left(\frac{2(n+2)}{\lambda(\lambda-1)R^{2}}\right)^n\int_{\Gamma_k}\frac{(\lambda R)^n M}{|t-z|}dt
  \\
  &\leq \sum_{k=0}^n \left(\frac{2\lambda(n+2)}{(\lambda-1)}\right)^{n+1}M=(n+1)\left(\frac{2\lambda(n+2)}{(\lambda-1)}\right)^{n+1}M
  \\
  &\leq \left(\frac{2\lambda(n+2)}{\lambda-1}\right)^{n+2}M,
  \end{align*}
  as $|z-t|\geq \frac{R(\lambda-1)}{2(n+2)}$, $|\Gamma_k|< 2\pi \lambda R$, and  $\lambda/(\lambda-1)\geq 1$.
  Equation \eqref{eq:est-F-minus} then follows from $|F(z)-\overline{z}^n F_n(z)|\leq M+|z^n F_n(z)|$, \eqref{eq:est-F_n} and the fact that $\lambda(n+2)/(\lambda-1)>1$.
   \pbox
   
 \noindent {\bf Proof of Lemma~\ref{max-mod-principle}
:}
  As $\widetilde{F}(z):=F(z)-F_n(z)\overline{z}^n$ is a polyanalytic function of order $n-1$ with 
  $$
  \sup_{\D(0,(1+\lambda)R/2)}|\widetilde{F}|\leq 2M\left(\frac{2\lambda(n+2)}{\lambda-1}\right)^{n+2},
  $$
  we can reapply Lemma~\ref{estimate-F_k} with $\lambda_1=(1+\lambda)/2$ to obtain
\begin{align*}
  |z^{n-1}F_{n-1}(z)|&= |z^{n-1}\widetilde{F}_{n-1}(z)|\leq 2M\left(\frac{2\lambda_1(n+1)}{\lambda_1-1}\right)^{n+1} \left(\frac{2\lambda(n+2)}{\lambda-1}\right)^{n+2},
 \end{align*}
and 
$$
|F(z)-\overline{z}^{n-1}F_{n-1}(z)-\overline{z}^nF_n(z)|\leq 4M\left(\frac{2\lambda_1(n+1)}{\lambda_1-1}\right)^{n+1} \left(\frac{2\lambda(n+2)}{\lambda-1}\right)^{n+2},
$$
for every $z\in\D(0,(1+\lambda_1)R/2)$.
Iterating this argument with 
$\lambda_l=\frac{1+\lambda_{l-1}}{2}=1+\frac{\lambda-1}{2^l}$ and setting $\lambda_0=\lambda$ then yields
\begin{align*}
  |z^{n-k}F_{n-k}(z)|&\leq 2^{k}M\prod_{l=0}^k\left(\frac{2\lambda_l(n+2-l)}{\lambda_l-1}\right)^{n+2-l}\leq  2^{n}M\prod_{l=0}^n\left(\frac{2\lambda_0(n+2-l)}{\lambda_n-1}\right)^{n+2}\\
  &=2^{n}M \left(\frac{2^{n+1}\lambda(n+2)!}{\lambda-1}\right)^{n+2}\leq M \left(\frac{2^{n+2}\lambda(n+2)!}{\lambda-1}\right)^{n+2}\\
 &\leq M \left(\frac{2\lambda}{\lambda-1}(n+2)^{n+2}\right)^{n+2} ,
 \end{align*}
 where we used that $\lambda_n-1=\frac{\lambda-1}{2^n}$ and that $2^{n-1}n!\leq n^n$. This is precisely the statement of Lemma~\ref{max-mod-principle}. \pbox

\subsection{Gabor Frames}
Let $\Gamma=\{z_i\}_{i\in\I}\subset \C$ be discrete.
A collection $\{\pi(z_i)g\}_{i\in\I}$ of time-frequency shifts of a window $g\in L^2(\R)$ is called  a 
\emph{Gabor frame},  if there exist constants $A,B>0$, called the frame bounds, such that
$$
A\|f\|^2\leq\sum_{i\in\I}|\langle f,\pi(z_i)g\rangle|^2\leq B\|f\|^2,\ \ \forall f\in L^2(\R).
$$
If only the right inequality is satisfied, then $\{\pi(z_i)g\}_{i\in\I}$ is called a \emph{Gabor Bessel sequence}.
A discrete set $\Gamma\subset\C$ is said to be \emph{uniformly separated} if $\inf\{|z-y|:\ z,y\in \Gamma,\ z\neq y\}>0$  and \emph{relatively uniformly separated} if it is the union of finitely many uniformly separated sets. 
The \emph{lower Beurling density} is defined as
$$
D^- (\Gamma):=\liminf_{R\rightarrow\infty}\inf_{z\in\C}\frac{\# \{\Gamma\cap z+ Q_R\}}{R^2}.
$$
The following result is due to Christensen, Deng and Heil
\cite{chdehe99}.
\begin{lemma}\label{eq:nec-density}
If $\Gamma$ generates  a Gabor frame, then $\Gamma$ is relatively uniformly separated and $D^- (\Gamma)\geq 1$.
\end{lemma}

%Throughout this contribution we will be particularly interested in irregular sets $\Gamma$ and connect the discrete frame property with lower bounds for continuous sub-sampling of the short-time Fourier transform.

%\begin{lemma}
%Let $0<a_1<a_2$. If $E$ is $(\gamma,a_1)$-dense, then $E$ is $\Big(\big(\frac{a_1}{a_2}\big)^2\gamma,a_2\Big)$-dense. In particular, if every $(\gamma,a_2)$-dense set $E$ is a planar set of sampling with asymptotic behavior  $ C_{2}\geq h(\gamma)$, then  every $(\gamma,a_1)$-dense set $E$ is a planar set of sampling with asymptotic behavior $C_{1}\geq h\big(\big(\frac{a_1}{a_2}\big)^2\gamma\big)$.
%\end{lemma}
%\proof First, observe that
%$$ \gamma=\inf_{(x,\omega)\in\R^2}\frac{|E\cap Q^{a_1}+(x,\omega)|}{a_1^2}\leq \left(\frac{a_2}{a_1}\right)^2\inf_{(x,\omega)\in\R^2}\frac{|E\cap Q^{a_2}+(x,\omega)|}{a_2^2}. $$
%The second part of the Lemma then follows directly from this observation as every $(\gamma,a_1)$-dense set is also    $\Big(\big(\frac{a_1}{a_2}\big)^2\gamma,a_2\Big)$-dense and thus a planar set of sampling.\pbox

\section{A Remez-Type Inequality for Polyanalytic Functions}

In this section, we derive a Remez-type inequality for polyanalytic functions. The proof relies on the the following result about plurisubharmonic functions, see
\cite[Theorem 1.2]{br99}.

\begin{theorem}[Brudnyi]\label{brudnyi}
Let $h:\C^d\rightarrow\R$ be plurisubharmonic, and let $r>1$. If $h$ satisfies 
$$
\sup_{B_\C(0,r)}h=0\quad\mbox{and}\quad
\sup_{B_\C(0,1)}h\geq -1,
$$
and $a>1$ is chosen such that $B_\R(x,t)\subset B_\C(x,at)\subset B_\C(0,1)$, then there exist constants $c=c(a,r)>0$ and $\kappa=\kappa(d)\geq 1$ such that the inequality 
\begin{equation}\label{remez-ineq-psh}
\sup_{B_\R(x,t)} h\leq c\ln\left(\frac{\kappa|B_\R(x,t)|}{|\Omega|}\right)+\sup_\Omega h,
\end{equation}
holds for every measurable set $\Omega\subset B_\R(x,t)$ of positive measure.
\end{theorem}
Note that, if a function $\Phi:\C^d\rightarrow\C$ is analytic, then $\ln|\Phi|$ is plurisubharmonic, see \cite[Section 2.2, p. 85]{kra82}.
Let $F$ be a polyanalytic function of order $n$ written as in \eqref{poly1}:
$$
F(z)=F(z,\overline{z})=\sum_{k=0}^n F_k(z)\overline{z}^k
$$
The function
\begin{equation}\label{f-to-F}
\Phi(F)(z_1,z_2):= \sum_{k=0}^n F_k(z_1+iz_2)(z_1-iz_2)^k,\qquad \forall (z_1,z_2)\in\C^2,
\end{equation}
is a holomorphic function in two complex variables. In particular, $\ln|\Phi(F)|$ is plurisubharmonic. If we set $z=x+iy$, $z_1=x$ and $z_2=y$ then $F(z)=F(x+iy)=\Phi(F)(x,y)=\Phi(F)(z_1,z_2)$, i.e. $F=\Phi(F)|_{\R\times\R}$.

Using Lemma~\ref{max-mod-principle}, we can estimate the supremum of $|\Phi(F)|$ in terms of a supremum of $|F|$.
\begin{lemma}\label{lem:first-estimate-of-poly}
If $F$ is a polyanalytic function of order $n$ with $M:=\sup_{\D(0,4R)}|F|$, 
then  
\begin{equation}\label{estimate-of-Phi}
\sup_{B_\C(0,2R)}|\Phi(F)(z_1,z_2)|
 \leq \left(4(n+2)\right)^{(n+2)^2}M.
\end{equation}
\end{lemma}
\proof First, note that if $(z_1,z_2)\in B_\C(0,2R)$, then $$|z_1\pm iz_2|^2\leq  (|z_1|+|z_2|)^2\leq 2(|z_1|^2+|z_2|^2)\leq 8R^2.$$ As $F_k$ is analytic, it follows by the maximum modulus principle that it attains its maximum at the boundary of the disc $\D(0,\sqrt{8}R)$ and consequently that
$$
 \sup_{\D(0,\sqrt{8}R)}|F_k(z)|(\sqrt{8}R)^k= \sup_{\D(0,\sqrt{8}R)}|F_k(z)z^k|\leq D_nM,
$$  
with Lemma~\ref{max-mod-principle}. Here, $D_n=D_n(4/\sqrt{8})=D_n(\sqrt{2})$. We  then conclude  that 
\begin{align*}
\sup_{B_\C(0,2R)}|\Phi(F)(z_1,z_2)|&\leq \sum_{k=0}^n\sup_{B_\C(0,2R)}|F_k(z_1+iz_2)|\sup_{B_\C(0,2R)}|z_1-iz_2|^k \nonumber\\
&\leq  \sum_{k=0}^n\sup_{\D(0,\sqrt{8}R)}|F_k(z)|(\sqrt{8}R)^k
\\
&\leq (n+1) D_n M\leq (n+1)8^{n+2}(n+2)^{(n+2)^2}M\\
&\leq  (4(n+2))^{(n+2)^2}M,
\end{align*}
where we used that $2\sqrt{2}/(\sqrt{2}-1)<8$ and $(n+1)8^{n+2}\leq 4^{n(n+2)}8^{n+2}=4^{(n+2)^2}$.
\pbox

\begin{theorem}\label{prop-poly-remez}
Let $0<\rho\leq R$, and $\Omega\subset \D(0,R)$ be measurable,
and  $F$ be polyanalytic of degree $n$ in $\D(0,5R)$.   If we write  $M:=\sup_{\D(0,4R)}|F|,$ and assume that $m:=|F(0)|>0$, then 
\begin{equation}
\sup_{\D(0,\rho)}|F|\leq \left(\frac{\kappa|\D(0,\rho)|}{|\Omega|}\right)^{c\left[\ln\frac{M}{m}+(n+2)^2\ln 4(n+2)\right]} \sup_{\Omega}|F|,
\end{equation}
where the constant $c$ can be chosen so that it only depends on the fraction $\rho/R$, and $\kappa$ is the constant in Brudnyi's theorem for $d=2$.
\end{theorem}
\proof Let $\alpha:=|\Phi(F)(0,0)|$ and $\beta:=\sup_{B_\C(0,2R)}|\Phi(F)|$. The function 
$$
h(z_1,z_2):=\frac{1}{\ln\frac{\beta}{\alpha}}\ln\frac{|\Phi(F)(\lambda Rz_1,\lambda Rz_2)|}{\beta},
$$
$1<\lambda <2$,
 is plurisubharmonic and satisfies 
$$
h(0,0)=\frac{1}{\ln\frac{\beta}{\alpha}}\ln\frac{|\Phi(F)(0,0)|}{\beta} \frac{\ln\frac{\alpha}{\beta}}{\ln\frac{\beta}{\alpha}}=-1,
$$
 and 
  $$
  \sup_{B_\C(0,\frac{2}{\lambda})}h=\frac{1}{\ln\frac{\beta}{\alpha}}\sup_{B_\C(0,2R)}\ln\frac{|\Phi(F)|}{\beta}=0.
  $$  
Hence, the assumptions of Theorem~\ref{brudnyi} are satisfied with $r=\frac{2}{\lambda}>1$, $d=2$ and $a=1/t$. Choosing  $t=\frac{\rho}{\lambda R}<1$, it therefore follows that there exist $c=c\big(\frac{\lambda R}{\rho},\frac{2}{\lambda}\big)$ and $\kappa=\kappa(2)$ such that 
\begin{align*}
\sup_{\D(0,\rho)} \frac{1}{\ln\frac{\beta}{\alpha}}\ln\frac{|F|}{\beta}&=\sup_{\D(0,\lambda Rt)} \frac{1}{\ln\frac{\beta}{\alpha}}\ln\frac{|F|}{\beta}=
\sup_{B_\R(0,\lambda Rt)} \frac{1}{\ln\frac{\beta}{\alpha}}\ln\frac{|\Phi(F)|}{\beta}
\\
&
=\sup_{B_\R(0,t)} \frac{1}{\ln\frac{\beta}{\alpha}}\ln\frac{|\Phi(F)(\lambda R\ \cdot\ )|}{\beta}
\\
&\leq c\ln\left(\frac{\kappa|B_\R(0,t)|}{|\Omega/\lambda R|}\right)+\sup_{\Omega/\lambda R} \frac{1}{\ln\frac{\beta}{\alpha}}\ln\frac{|\Phi(F)(\lambda R\ \cdot\ )|}{\beta}
\\
&= c\ln\left(\frac{\kappa \lambda^2R^2|\D(0,t)|}{|\Omega|}\right)+\sup_{\Omega} \frac{1}{\ln\frac{\beta}{\alpha}}\ln\frac{|\Phi(F)|}{\beta}
\\
&= c\ln\left(\frac{\kappa|\D(0,\lambda Rt)|}{|\Omega|}\right)+\sup_{\Omega} \frac{1}{\ln\frac{\beta}{\alpha}}\ln\frac{|F|}{\beta}.
\end{align*}
Taking the exponential function of both sides of the inequality and recalling the definition of $t$ yields
$$
\left(\sup_{\D(0,\rho)} \frac{|F|}{\beta}\right)^{\frac{1}{\ln\frac{\beta}{\alpha}}}\leq  \left(\frac{\kappa|\D(0,\rho)|}{|\Omega|}\right)^c \left(\sup_{\Omega} \frac{|F|}{\beta}\right)^{\frac{1}{\ln\frac{\beta}{\alpha}}},
$$
and consequently
$$
\sup_{\D(0,\rho)} |F|\leq  \left(\frac{\kappa|\D(0,\rho)|}{|\Omega|}\right)^{c\ln\frac{\beta}{\alpha}} \sup_{\Omega} |F|.
$$
This concludes the proof when  observing that $\alpha=m$ and that $\beta\leq (4(n+2))^{(n+2)^2}M$ by Lemma~\ref{lem:first-estimate-of-poly}.
\pbox

\section{Lower Sampling Bounds for the STFT with Hermite Windows}

This section is devoted to the proof of Theorem~\ref{thm:main}. 
Before going into the details, let us shortly state a direct consequence of this result for functions in the true polyanalytic Bargmann spaces $\F_p^n$.

\begin{corollary}\label{cor:main}
Let $1\leq p<\infty$, and $F\in \F^n_p(\C)$.
If $\Omega\subset \C$ is $(\gamma,R)$-dense for some scale $R>0$, then there exists $\eta=\eta(n,R)$ and $\sigma=\sigma(n,R)$ and a numerical constant $C>0$ such that
\begin{equation}
\|F\|_{\mathcal{L}^p(\C)}\leq \eta\left(\frac{\gamma}{C}\right)^{-  \sigma}\|F\|_{\mathcal{L}^p(\Omega)}.
\end{equation}
\end{corollary}
\proof The result follows from Theorem~\ref{thm:main} once we recall that $|\B^{n+1}f(z)|e^{-\pi|z|^2/2}=|V_{h_n}f(\overline{z})|$, $\B^{n+1}(M^p(\R))=\F^n_p(\C)$, and the definition of $\mathcal{L}^p(\Omega)$ in \eqref{def:LpC}.\pbox

\noindent Note that  if we set $\rho=R$, and make the particular choice $\lambda=\sqrt{2}$ in Theorem~\ref{prop-poly-remez}, it follows that the constant  $c=c(\sqrt{2},\sqrt{2})$  is independent of $R$. We are now in place to show a Remez type inequality for the short-time Fourier transform with Hermite windows.

\begin{lemma}\label{lem-infty-rem-stft}
Let $w\in\C$, $\Omega\subset \D(w,R)$ be measurable, $m:=|V_{h_n}f(w)|>0$, and $$\gamma:=\int_{\D(w,5R)} |V_{h_n}f(z)| dz.$$ There exists numerical constants $c>0$ and $\kappa\geq 1$ (independent of $R$, $n$ and $f$) such that, for every $f\in M^\infty(\R)$,
\begin{equation}
\sup_{z\in \D(w,R)}|V_{h_n}f(z)|\leq e^{\frac{\pi}{2} R^2}\left(\frac{\kappa|\D(0,R)|}{|\Omega|}\right)^{K(R,n,\gamma,m)}\sup_{z\in \Omega}|V_{h_n}f(z)|,
\end{equation}
where $K(R,n,\gamma,m):=c\left[\ln\frac{\gamma}{m}+8\pi R^2+\ln( \nu_n(R)^{-1})+(n+2)^2\ln 4(n+2)\right]$
and $\nu_n(R)$ is given by \eqref{eq:local-const}.
\end{lemma}
\proof As $|V_{h_n}f(z-w)|=|V_{h_n}\pi(w)f(z)|$, we may without loss of generality assume that $w=0$. Set $M:=\sup_{\D(0,4R)}|\B^{n+1}f|$. From \eqref{def-poly-bargmann-trafo} we know that $m=|V_{h_n}f(0)|=|\mathcal{B}^{n+1}f(0)|$ and that 
\begin{equation}\label{eq:abs-berg-stft}
|V_{h_n}f(z)|=|\mathcal{B}^{n+1}f(\overline{z})|e^{-\pi|z|^2/2}.
\end{equation}
 Let us first estimate $M$ in terms of the quantity $\gamma$.
By the local reproducing formula \eqref{eq:local-repr}, we have that
\begin{align}\label{eq:est-M}
M&=\sup_{z\in\D(0,4R)}|\B^{n+1}f(z)|=\sup_{z\in \D(0,4R)}|V_{h_n}f(z)e^{\pi |z|^2/2}|\nonumber
\\
&\leq e^{8\pi R^2}\sup_{z\in \D(0,4R)} \nu_n(R)^{-1}\int_{\D(z,R)} |V_{h_n}f(w)\langle \pi(w)h_n,\pi(z)h_n\rangle| dw\nonumber
\\
&\leq e^{8\pi R^2} \nu_n(R)^{-1} \int_{\D(0,5R)} |V_{h_n}f(w)| dw=e^{8\pi R^2} \nu_n(R)^{-1}\gamma.
\end{align}
Let $\Omega^\ast:=\{\overline{z}:\ z\in\Omega\}$. As $\B^{n+1}f$ is polyanalytic of order $n$, we may use Theorem~\ref{prop-poly-remez} and \eqref{eq:abs-berg-stft} to show that
\begin{align*}
\sup_{z\in \D(0,R)}|V_{h_n}f(z)|&=\sup_{z\in \D(0,R)}|\mathcal{B}^{n+1}f(z)|e^{-\pi|z|^2/2}\leq \sup_{z\in \D(0,R)}|\mathcal{B}^{n+1}f(z)|
\\
&\leq \left(\frac{\kappa|\D(0,R)|}{|\Omega^\ast|}\right)^{c\left[\ln\frac{M}{m}+(n+2)^2\ln 4(n+2)\right]}\sup_{z\in \Omega^\ast}|\mathcal{B}^{n+1}f(z)|
\\
&\leq e^{\frac{\pi}{2} R^2}\left(\frac{\kappa|\D(0,R)|}{|\Omega|}\right)^{c\left[\ln\frac{M}{m}+(n+2)^2\ln 4(n+2)\right]}\sup_{z\in \Omega^\ast}\big|\mathcal{B}^{n+1}f(z)\big|e^{-\pi|z|^2/2}.
\end{align*}
Using again \eqref{eq:abs-berg-stft}, we obtain $\sup_{z\in \Omega^\ast}|\mathcal{B}^{n+1}f(z)|e^{-\pi|z|^2/2}=\sup_{z\in \Omega}|V_{h_n}f(z)|$ and the result follows once we plug in the estimate from \eqref{eq:est-M}. \pbox

\begin{lemma}\label{Lp-estimate}
Let $1\leq p<\infty$.
With the notation and conditions of Lemma~\ref{lem-infty-rem-stft}
$$
 \|V_{h_n} f\|_{L^p(\D(w,R))}\leq e^{ \pi R^2/2}\left(\frac{2\kappa|\D(0,R)|}{|\Omega|}\right)^{K(R,n,\gamma,m)+1}\|V_{h_n} f\|_{L^p(\Omega)}
$$
for every $f\in M^\infty(\R)$.
\end{lemma}
\proof As before, we may assume that $w=0$. For $\theta>0$ we define the set
$$
A_\theta:=\left\{w\in \D(0,R):\ |V_{h_n} f(w)|<e^{-\pi R^2/2}\left(\frac{\theta}{\kappa|\D(0,R)|}\right)^{K(R,n,\gamma,m)}\sup_{\D(0,R)}|V_\varphi f|\right\}.
$$
Taking $\Omega=A_\theta$ in  Lemma~\ref{lem-infty-rem-stft} yields 
\begin{align*}
\left(\frac{|A_\theta|}{\kappa|\D(0,R)|}\right)^{K(R,n,\gamma,m)}\sup_{\D(0,R)}|V_{h_n} f|&\leq e^{\pi R^2/2}\sup_{ A_\theta}|V_{h_n} f|\\ &<\left(\frac{\theta}{\kappa|\D(0,R)|}\right)^{K(R,n,\gamma,m)}\sup_{\D(0,R)}|V_{h_n}f|,
\end{align*}
which shows that $|A_\theta|< \theta$.
Let now $\Omega\subset \D(0,R)$, and $\theta=\frac{|\Omega|}{2}$. Then $\frac{|\Omega|}{2}\leq |\Omega\cap A_{|\Omega|/2}^c|$ and, as $\kappa\geq 1$, we obtain
\begin{align*}
e^{ -\pi pR^2/2}&\left(\frac{|\Omega|}{2\kappa|\D(0,R)|}\right)^{pK(R,n,\gamma)+p}\int_{\D(0,R)}|V_{h_n} f(w)|^pdw\\
&\leq \kappa e^{ - \pi pR^2/2}\left(\frac{|\Omega|}{2\kappa|\D(0,R)|}\right)^{p\cdot K(R,n,\gamma,m)+1}\int_{\D(0,R)}|V_{h_n} f(w)|^pdw
\\
&\leq \frac{|\Omega|}{2}e^{ -\pi pR^2/2}\left(\frac{|\Omega|}{2\kappa|\D(0,R)|}\right)^{p\cdot K(R,n,\gamma,m)}\sup_{\D(0,R)}|V_{h_n} f|^p
\\
&\leq \int_{\Omega\cap A_{|\Omega/2|}^c}e^{ -\pi pR^2/2}\left(\frac{|\Omega|}{2\kappa|\D(0,R)|}\right)^{p\cdot K(R,n,\gamma,m)}\sup_{\D(0,R)}|V_{h_n} f|^pdw
\\
&\leq \int_{\Omega \cap A_{|\Omega|/2}^c} |V_{h_n} f(w)|^pdw\leq \int_\Omega |V_{h_n}f(w)|^pdw,
\end{align*}
where we used the definition of $A_{|\Omega/2|}^c$ to derive the second to last inequality. 
\pbox

\begin{lemma}\label{size-of-W-theta}
Let $1\leq p<\infty$, and $\theta,R>0$. If  the set $W_{\theta,R}\subset\C$ is defined as 
$$
W_{\theta,R}:=\left\{z\in \C:\ |V_{h_n} f(z)|^p\leq \frac{\theta}{|\D(0,R)|}\int_{\D(z,R)}|V_{h_n} f(w)|^p dw\right\},
$$
then the following inequality holds
\begin{equation}\label{eq-W_e}
\int_{W_{\theta,R}^c}|V_{h_n} f(z)|^p dz\geq(1-\theta)\int_{\C}|V_{h_n} f(z)|^pdz,\ \ \forall f\in M^p(\R).
\end{equation}
\end{lemma}
\proof
By definition of $W_{\theta,R}$ and Fubini's theorem we can write
\begin{align*}
 \int_{W_{\theta,R}}|V_{h_n} f(z)|^pdz&\leq \int_{W_{\theta,R}}\frac{\theta}{|\D(0,R)|}\int_{\D(z,R)}|V_{h_n} f(w)|^p dwdz
\\
&=
\frac{\theta}{|\D(0,R)|} \int_{\C}|V_{h_n} f(w)|^p\int_{W_{\theta,R}}\chi_{\D(w,R)}(z)dz dw 
\\
&\leq \theta \int_{\C}|V_{h_n} f(w)|^pdw,
\end{align*}
as claimed.
 \pbox

 \begin{lemma}\label{lem:W-theta-estimates-beta}
 If $1\leq p<\infty$ and $w\in W_{1/2^p,5R}^c$, then there exists a constant ${C}>0$ such that 
 $$ K(R,n,\gamma,m)\leq  C\left(R^2 +\ln( \nu_n(R)^{-1}) +n^2\ln n +1\right)=:b(R,n)-1,$$ which is independent of $f$.
 \end{lemma}
 \proof It is enough to show that $\gamma/m$ is bounded independent of $f$. By the definition of $ W_{\theta,5R}^c$ it follows by H{\"o}lder's inequality that if $w\in  W_{1/2^p,5R}^c$ we have
\begin{align*}
\gamma^p&=\left(\int_{\D(w,5R)}|V_{h_n}f(z)|dz\right)^p\leq |\D(0,5R)|^{p/q}\int_{\D(w,5R)}|V_{h_n}f(z)|^pdz
\\
&< \frac{|\D(0,5R)|^{1+p/q}}{1/2^p}|V_{h_n}f(w)|^p=2^p|\D(0,5R)|^{p}m^p.
\end{align*}
 Consequently, $\frac{\gamma}{m}< 2|\D(0,5R)|=50\pi R^2,$ and 
\begin{align*}
K(R,n,\gamma,m)&\leq c\left[ \ln 50\pi R^2 +8\pi R^2 +\ln (\nu_n(R)^{-1})+(n+2)^2\ln 4(n+2) \right]
\\
&\leq c\left( \widetilde C R^2 +\ln(\nu_n(R)^{-1}) +(n+2)^2\ln 4(n+2)\right)
\\
&\leq C\left(  R^2 +\ln(\nu_n(R)^{-1}) +n^2\ln n+1\right),
\end{align*}
as claimed.
 \pbox
 
\noindent We are now in position to prove Theorem~\ref{thm:main} which we restate here in a more precise form:

\begin{theorem}\label{thm:mainprecise}
Let $1\leq p<\infty$, $f\in M^p(\R)$.
If $\Omega\subset \C$ is $(\gamma,R)$-dense for some scale $R>0$, then there exists a numerical constant $C>0$  such that 
\begin{equation}
\|V_{h_n}f\|_{L^p(\C)}\leq \eta(n,R)\left(\frac{\gamma}{C}\right)^{-  \sigma(R,n)}\|V_{h_n}f\|_{L^p(\Omega)},
\end{equation}
where 
$$\sigma(n,R):=C(R^2+\ln(\nu_n(R)^{-1})+n^2\ln n+1),\ \  \mbox{ 
and } \ \ \ 
\eta(n,R):= \frac{R^2}{\nu_n(R)}C^{R^2+1}.
$$
\end{theorem}
\proof By Fubini's theorem and  Lemma~\ref{Lp-estimate} we may derive
\begin{align*} 
\int_\Omega|V_{h_n} f(z)|^pdz& \geq \int_\Omega|V_{h_n} f(z)|^p\frac{1}{|\D(0,R)|}\int_{W_{1/2^p,5R}^c}\chi_{\D(z,R)}(w)dwdz
\\
&=\int_{W_{1/2^p,5R}^c} \frac{1}{|\D(0,R)|}\int_{\Omega\cap \D(w,R)}|V_{h_n}  f(z)|^p dzdw
\\
&\geq \frac{e^{-\pi pR^2/2}}{|\D(0,R)|} \int_{W_{1/2^p,5R}^c}\hspace{-0.15cm} \left(\frac{|\Omega\cap \D(w,R)}{ 2\kappa|\D(w,R)|}\right)^{p(K(R,n,\gamma,m)+1)}\hspace{-0.2cm} \int_{\D(w,R)}|V_{h_n}  f(z)|^p dzdw.
\end{align*}
Now,  Lemma~\ref{lem:W-theta-estimates-beta} allows to estimate $K(R,n,\gamma,m)$. As $b(R,n)$ is independent of $f$ and $w$, and $\Omega$ is $(\gamma,R)$-dense, it follows that 
\begin{equation}\label{eq:forelast-eq}
\int_\Omega|V_{h_n} f(z)|^pdz \geq \frac{e^{- \pi p R^2/2}}{|\D(0,R)|} \left(\frac{\gamma}{2\kappa }\right)^{p\cdot b(R,n)}\int_{W_{1/2^p,5R}^c} \int_{\D(w,R)}|V_{h_n}  f(z)|^p dzdw.
\end{equation}
It remains to estimate the double integral on the right hand side. H{\"o}lder's inequality, the local reproducing formula \eqref{eq:local-repr}, and Lemma \ref{size-of-W-theta} give
\begin{align*}
\int_{W_{1/2^p,5R}^c} \int_{\D(w,R)}&|V_{h_n}  f(z)|^p dzdw \geq
\int_{W_{1/2^p,5R}^c} \int_{\D(w,R)}|V_{h_n}  f(z)\langle \pi(z)h_n,\pi(w)  h_n\rangle|^p dzdw
\\
 &\geq \int_{W_{1/2^p,5R}^c}\hspace{-0.4cm} |\D(0,R)|^{-p/q}\left|\int_{\D(w,R)}\hspace{-0.3cm}V_{h_n}  f(z)\langle \pi(z)h_n,\pi(w)  h_n\rangle dz\right|^p dw 
\\
&=  |\D(0,R)|^{-p/q}\int_{W_{1/2^p,5R}^c} \nu_n(R)^p|V_{h_n}  f(w)|^p dw 
\\
&\geq |\D(0,R)|^{-p/q}\nu_n(R)^p\left(1-2^{-p}\right)\int_{\C}|V_{h_n}  f(z)|^pdz 
\\
&\geq \left(\frac{\nu_n(R)}{2|\D(0,R)|^{p/q}}\right)^p \int_{\C}|V_{h_n}  f(z)|^pdz, 
\end{align*}
where we used $1-2^{-p}\geq 2^{-p}$ in the final step. Plugging this into \eqref{eq:forelast-eq} finally yields
$$
\int_\Omega|V_{h_n} f(z)|^pdz \geq \left(\frac{e^{-\pi R^2/2}\nu_n(R)}{2|\D(0,R)|}\right)^p \left(\frac{\gamma}{2\kappa }\right)^{p\cdot b(R,n)}\int_{\C}|V_{h_n}  f(z)|^pdz
$$
as claimed.\pbox

\section{Sufficient Density Conditions for General Windows}

\subsection{Irregular Gabor Frames}\label{sec:irreg-gabor}

There is only little literature on irregular Gabor frames. Coorbit theory \cite{fegr89,groe91}, on the one hand,  guarantees the frame property for "nice" windows and irregular sampling points with a sufficiently high density. However, the results do not provide any estimate of how to choose the density and how the frame bounds behave. 
Gr{\"o}chenig on the other hand derived quantitative results in \cite{groe93}. The choice of the sampling sets  however do not leave enough freedom for our purposes in this section. It is a result by Sun and Zhou  \cite[Lemma 2.6]{suzho02} that has both ingredients: quantitative estimates of the frame bounds and enough freedom in choosing the sampling points.
Recall that the
Sobolev spaces $H^s(\R)$ are defined as
$$
H^s(\R):=\left\{f\in L^2(\R):\ \int_\R (1+|\xi|^2)^s |\widehat{f}(\xi)|^2d\xi<\infty\right\}.
$$
\begin{theorem}[Sun \& Zhou]\label{lem:suzhou}
Let $g,tg\in H^1(\R)$ and let $R>0$ be such that 
\begin{equation}\label{eq:con-frame-suzhou}
\Delta:=\frac{2R}{\pi}\left(\|g'\|_2+\|tg\|_2+\frac{2R}{\pi}\|tg'\|_2\right)<\|g\|_2.
\end{equation}
Moreover, let $Q_n$ be a collection of squares with side length $R_n\leq R$ such that $\bigcup_{n\in\N }Q_n=\C$ and $|Q_n\cap Q_m|=0$, $n\neq m$. Then for any $z_n\in Q_n$, $\{R_n \pi(z_n)g\}_{n\in\N}$ is a frame for $L^2(\R)$ with frame bounds $A\geq\big(\|g\|_2-\Delta\big)^2$ and $B\leq \big(\|g\|_2+\Delta\big)^2$.
\end{theorem}
Theorem 2.8 in \cite{suzho02} gives a more detailed picture of the frame structure for this class of  windows. The full generality of the result is however not needed for our purposes.

\begin{corollary}\label{cor:frame-cs}
If $g\in H^1(\R)$ is compactly supported in $[-S,S]$, and $|x_{n,m}-Rn|<R/2$, $|\xi_{n,m}-Rm|<R/2$ and $R<\min\left(\frac{\pi \|g\|_2}{4\|g'\|_2},\frac{1}{2S}\right)$, then $\{\pi(z_{n,m})g\}_{n,m\in\Z}$ is a frame for $L^2(\R)$ with 
\begin{equation}\label{eq:lower-irr-cs}
\frac{1}{3 R^2}\left(\|g\|_2-\frac{4R}{\pi}\|g'\|_2\right)^2\|f\|^2_2\leq \sum_{n,m\in\Z}|V_gf(z_{n,m})|^2\leq \frac{2}{R^2}\left(\|g\|_2+\frac{2R}{\pi}\|g'\|_2\right)^2\|f\|_2^2.
\end{equation}
\end{corollary}
\proof First note that $g\in H^1(\R)$ compactly supported implies that $tg\in H^1(\R)$.
Using the assumption $R<1/2S$, $\Delta$ can be estimated as
\begin{align*}
\Delta&=\frac{2R}{\pi}\left(\|g'\|_2+\|tg\|_2+\frac{2R}{\pi}\|tg'\|_2\right)\leq \frac{2R}{\pi}\left(\|g'\|_2+S\|g\|_2+\frac{1}{\pi}\|g'\|_2\right)\\ 
&\leq \frac{2R}{\pi}\frac{\pi+ 1}{\pi}\|g'\|_2+\frac{1}{\pi}\|g\|_2.
\end{align*}
Hence, $\Delta<\|g\|_2$ if
%$$ \frac{2a}{\pi}\frac{\pi+ 1}{\pi}\|g'\|+\frac{1}{\pi}\|g\|<\|g\|,$$
%or equivalently 
$$
R<\frac{\pi\|g\|_2}{4\|g'\|_2}<\frac{\pi(\pi- 1)\|g\|_2}{2(\pi+1)\|g'\|_2},
$$
as $1/2<\frac{\pi-1}{\pi+1}$. % sets the sufficient condition on $a$. 
Using \eqref{eq:con-frame-suzhou} and the estimate for $\Delta$, the upper frame bound is estimated by
$$
B\leq \frac{1}{R^2}\left(\|g\|_2+\Delta\right)^2\leq \frac{1}{R^2}\left(\frac{\pi+1}{\pi}\right)^2\left(\|g\|_2+\frac{2R}{\pi}\|g'\|_2\right)^2\leq \frac{2}{R^2}\left(\|g\|_2+\frac{2R}{\pi}\|g'\|_2\right)^2.
$$
For the lower frame bound, we have
\begin{align*}
A&\geq \frac{1}{R^2}\left(\|g\|_2-\Delta\right)^2
=\frac{1}{R^2}\left(\frac{\pi-1}{\pi}\right)^2\left(\|g\|_2-\frac{2R}{\pi}\frac{\pi+1}{\pi-1}\|g'\|_2\right)^2
\\
&\geq \frac{1}{3R^2}\left(\|g\|_2-\frac{4R}{\pi}\|g'\|_2\right)^2
\end{align*}
as claimed.
\pbox

\subsection{Sufficient Density Conditions for Planar Sets of Sampling}\label{sec:planar-sampling}

For the rest of this paper let us slightly change the definition of $(\gamma,R)$-dense sets. Instead of discs we now use squares to define such sets. In particular,
a set $\Omega\subset \C$ is called 
\emph{$(\gamma,R)$-dense} if
\begin{equation}\label{eq:gR-dense-Q}
\gamma:=\inf_{z\in\C}\frac{|\Omega\cap z+Q_R|}{R^2}>0.
\end{equation}
Our main goal in this section is to determine under which conditions a $(\gamma,R)$-dense set at  small enough scales $R$ is a  planar set of sampling. It is therefore irrelevant which definition of density we use as every $(\gamma,R)$-dense set in the sense of \eqref{eq:gR-dense-D} is $(\gamma/2,R)$-dense set in the sense of \eqref{eq:gR-dense-Q} and every $(\gamma,R)$-dense set in the sense of \eqref{eq:gR-dense-Q} is $(\gamma/2,\sqrt{2}R)$-dense set in the sense of \eqref{eq:gR-dense-D}.

Let us write $Q_R({n,m}):=Rn+iRm+Q_R$.
We are now able to establish a connection between irregular Gabor frames and  planar sets of sampling.

\begin{theorem}\label{thm:planar-sampling-frames}
Let $g\in L^2(\R)$, and $\Omega\subset \C$ be measurable. 
\begin{enumerate}[(i)]
\item If there exists $R>0$, such that $\{\pi({z_{n,m}})g\}_{n,m\in\Z}$ is a frame for every choice $z_{n,m}\in Q_R({n,m})$  with a lower frame bound $A$ independent of the particular choice of sampling points,
then every $(\gamma,R)$-dense set   $\Omega$ is a planar set of sampling. The sampling bound satisfies $C\leq \frac{\|g\|^2_2}{AR^2}\gamma^{-1}$.
\end{enumerate}
If in addition, $g$ is such that every relatively uniformly separated sequence of points in $\C$ generates a Gabor Bessel sequence, then the following statements hold true. 
\begin{enumerate}[(i)]
\setcounter{enumi}{1}
\item If there exists $R>0$ such that every $(\gamma,R)$-dense set $\Omega$ is a planar set of sampling with sampling bound $C\leq K \gamma^{-1}$, then $\{\pi(z_{n,m})g\}_{n,m\in\Z}$ is a frame for every choice $z_{n,m}\in Q_{R/2}({n,m})$ with lower frame bound $A\geq\frac{\|g\|^2_2}{KR^2}$.
\item There is no scale $R>0$ and $\varepsilon >0$ such that the sampling bound satisfies $C\leq K \gamma^{-1+\varepsilon}$ for every $(\gamma,R)$-dense set $\Omega$.
\end{enumerate}
\end{theorem}
\proof \emph{Ad $(i)$:} Let $f\in L^2(\R)$ be fixed and let $\Omega$ be a $(\gamma,R)$-dense set. For every $n,m\in\Z$ there exist $z_{n,m}\in \Omega\cap Q_R({n,m})$ such that 
$$
|V_gf(z_{n,m})|^2\leq \frac{1}{|\Omega\cap Q_R({n,m})|}\int_{\Omega\cap Q_R({n,m})}|V_gf(z)|^2dz.
$$
If $A$ is the uniform lower frame bound, then
\begin{align*}
\|V_gf\|^2_{L^2(\C)}&=\|g\|^2_2\|f\|^2_2\leq \frac{\|g\|_2^2}{A}  \sum_{n,m\in\Z}|V_gf(z_{n,m})|^2
\\
&\leq \frac{\|g\|^2_2}{AR^2}\sum_{n,m\in\Z} \frac{R^2}{|\Omega\cap Q_R({n,m})|} \int_{\Omega\cap Q_R({n,m})}|V_gf(z)|^2dz
\\
&\leq \frac{\|g\|^2_2}{AR^2}\gamma^{-1}\sum_{n,m\in\Z} \int_{\Omega\cap Q_R({n,m})}|V_gf(z)|^2dz
\\
&= \frac{\|g\|^2_2}{AR^2} \gamma^{-1} \int_{\Omega}|V_gf(z)|^2dz.
\end{align*}
\emph{Ad $(ii)$:}  Let $z_{n,m}\in Q_{R/2}({n,m})$ and choose $\Omega_\mu$ such that $\Omega_\mu\cap Q_{R/2}({n,m})$ is an  open neighborhood of $z_{n,m}$ and that $|\Omega_\mu\cap Q_{R/2}({n,m})|=\mu R^2$. It then follows that $\Omega_\mu$ is a  $(\gamma,R)$-dense set with $\mu\leq\gamma$, since every square $z+Q_R$ contains at least one square $Q_{R/2}({n,m})$. As $C\leq K \gamma^{-1}$, we have
$$
\int_{\C}|V_gf(z)|^2dz\leq K \gamma^{-1}\int_{\Omega_\mu}|V_gf(z)|^2dz \leq K \mu^{-1}\int_{\Omega_\mu}|V_gf(z)|^2dz
$$
$$
= K\sum_{n,m\in\Z}\mu^{-1}\int_{\Omega_\mu\cap Q_{R/2}({n,m})}|V_gf(z)|^2dz.
$$
Moreover, since  
\begin{align*}
\frac{1}{ \mu}\int_{\Omega_\mu\cap Q_{R/2}({n,m})}|V_gf(z)|^2dz&=\frac{R^2}{ |\Omega_\mu\cap Q_{R/2}({n,m})|}\int_{\Omega_\mu\cap Q_{R/2}({n,m})}|V_gf(z)|^2dz
\\
&\leq R^2 \sup_{z\in  Q_{R/2}({n,m})}|V_gf(z)|^2,
\end{align*} 
it follows by the assumption on $g$ that, when considering the limit $\mu\rightarrow 0$, we may apply the dominated convergence theorem to obtain
\begin{align*}
\|g\|^2_2\|f\|^2_2&=\int_{\C}|V_gf(z)|^2dz\leq K\lim_{\mu\rightarrow 0}\sum_{n,m\in\Z}\frac{1}{\mu}\int_{\Omega_\mu\cap Q_{R/2}({n,m})}|V_gf(z)|^2dz
\\ 
&=K\sum_{n,m\in\Z}\lim_{\mu\rightarrow 0}\frac{R^2}{|\Omega_\mu\cap Q_{R/2}({n,m})|}\int_{\Omega_\mu\cap Q_{R/2}({n,m})}|V_gf(z)|^2dz
\\
&=R^2 K\sum_{n,m\in\Z}|V_gf(z_{n,m})|^2.
\end{align*}
\emph{Ad $(iii)$:} Assume to the contrary, that there exist $R,\varepsilon>0$ such that $C\leq K \gamma^{-1+\varepsilon}\leq K\gamma^{-1}$. In particular, the assumption of $(ii)$ is satisfied and an arbitrary choice $z_{n,m}\in Q_{R/2}({n,m})$  generates a Gabor frame. On the other hand, repeating the calculations of the proof of $(ii)$ with $K \gamma^{1-\varepsilon}$ instead of $K \gamma$ gives 
$$
\frac{\|g\|_2^2 }{KR^2}\gamma^{-\varepsilon}\|f\|^2\leq \sum_{n,m}|V_gf(z_{n,m})|^2,
$$
for any $\gamma \in(0,1]$. Taking the limit $\gamma\rightarrow 0$ then shows, that $z_{n,m}\in Q_{R/2}({n,m})$ cannot generate a Gabor frame, a contradiction.
\pbox

\noindent The necessary density condition for Gabor frames in Lemma~\ref{eq:nec-density} now  shows that the sampling bound $C$  be bounded by $\gamma^{-1}$ only  for small scales $R$.

\begin{corollary}
If there exists $R>2$ such that every $(\gamma,R)$-dense set $\Omega$ is a planar set of sampling, then there exists no constant $K>0$ such that the  sampling constant is bounded by $K\gamma^{-1}$. In particular, $C\geq K\gamma^{-1-\varepsilon}$ for some $\varepsilon>0$.
\end{corollary}
\noindent \proof
Assume to the contrary, that there exists $K>0$ such that the sampling bound satisfies $C\leq K\gamma^{-1}$. Then, by Theorem~\ref{thm:planar-sampling-frames}~$(ii)$, it follows that $\{M_{Rm/2}T_{{Rn/2}}g\}_{n,m\in\Z}$ is a frame. The set $\Gamma=\{Rn/2+iRm/2\}_{n,m\in\Z}$ set has lower Beurling density $D^-(\Gamma)=4/R^2$. Now since $R>2$, it follows that  $D^-(\Gamma)<1$ which, by Lemma~\ref{eq:nec-density}, contradicts the assumption that $\{M_{Rm/2}T_{Rn/2}g\}_{n,m\in\Z}$ is a frame.
\pbox

\noindent We can now prove Theorem \ref{thm:main2} from the introduction in a more precise form:

\begin{corollary}
Let $g\in H^1(\R)$ be supported in $[-S,S]$. If  
$$
R<\min\left(\frac{\pi\|g\|_2}{4\|g'\|_2},\frac{1}{2S}\right),
$$
then every $(\gamma,R)$-dense set $\Omega\subset \C$ is a planar set of sampling. In particular, for every $f\in L^2(\R)$, one has
\begin{equation}
  \int_{\C}|V_gf(z)|^2dz\leq \frac{3}{\gamma}\left(1-\frac{4\|g'\|_2}{\pi\|g\|_2}R\right)^{-2} \int_{\Omega}|V_gf(z)|^2dz.
\end{equation}
\end{corollary}
\proof This result follows directly from Theorem~\ref{thm:planar-sampling-frames}~$(i)$ once one observes that the sampling constant is less than $\frac{\|g\|^2_2}{AR^2}\gamma^{-1}$, where $A$ is the lower frame bound given in \eqref{eq:lower-irr-cs}. 
\pbox

\section*{Acknowledgements} 
M. Speckbacher was supported by an Erwin-Schr{\"{o}}dinger Fellowship (J-4254) of the Austrian Science Fund FWF.

\end{document}